\documentclass[preprint,12pt]{elsarticle}

\usepackage{amssymb}
\usepackage{amsmath,amssymb,amsfonts}
\usepackage{tikz}
\usepackage{tabularx}

\newtheorem{theorem}{Theorem}
\newtheorem{corollary}{Corollary}

\newtheorem{proposition}{Proposition}
\newtheorem{remark}{Remark}
\newproof{proof}{Proof}

\begin{document}

\begin{frontmatter}

\title{Exponential stability of damped Euler-Bernoulli beam controlled by boundary springs and dampers}

%
%

\author[1]{Onur Baysal\corref{cor3}\fnref{fn1}}
\ead{onur.baysal@um.edu.mt}
\author[2]{Alemdar Hasanov%
\fnref{fn2}}
\ead{alemdar.hasanoglu@gmail.com}
\author[3]{Alexandre Kawano\fnref{fn3}}
\ead{akawano@usp.br}

\cortext[cor3]{Corresponding author}
\fntext[fn1]{Department of Mathematics, University of Malta, Malta}
\fntext[fn2]{Department of Mathematics, Kocaeli University, Turkey}
\fntext[fn3]{Escola Polit\'{e}cnica, University of S\~{a}o Paulo, S\~{a}o Paulo 05508900, Brazil}

\address[1]{Department of Mathematics, University of Malta, Malta}
\address[2]{Kocaeli University, 41001, Kocaeli, Turkey}
\address[3]{University of S\~{a}o Paulo, S\~{a}o Paulo 05508900, Brazil}

\begin{abstract}
In this paper, the vibration model of an elastic beam, governed by the damped Euler-Bernoulli equation $\rho(x)u_{tt}+\mu(x)u_{t}$$+\left(r(x)u_{xx}\right)_{xx}=0$, subject to the clamped boundary conditions $u(0,t)=u_x(0,t)=0$ at $x=0$, and the boundary conditions $\left(-r(x)u_{xx}\right)_{x=\ell}=k_r u_x(\ell,t)+k_a u_{xt}(\ell,t)$,  $\left(-\left(r(x)u_{xx}\right)_{x}\right )_{x=\ell}$$=- k_d u(\ell,t)-k_v u_{t}(\ell,t)$ at $x=\ell$, is analyzed. The boundary conditions at $x=\ell$ correspond to linear combinations of damping moments caused by rotation and angular velocity and also, of forces caused by displacement and velocity, respectively. The system stability analysis based on well-known Lyapunov approach is developed. Under the natural assumptions guaranteeing the existence of a regular weak solution, uniform exponential decay estimate for the energy of the system is derived. The decay rate constant in this estimate depends only on the physical and geometric parameters of the beam, including the viscous external damping coefficient $\mu(x) \ge 0$, and the boundary springs $k_r,k_d \ge 0$ and dampers $k_a,k_v \ge 0$. Some numerical examples are given to illustrate the role of the damping coefficient and the boundary dampers.
\end{abstract}

\end{frontmatter}


\section{Introduction}
In his paper, we study the exponential stability of the system governed by the following initial boundary value problem for the non-homogeneous damped Euler-Bernoulli beam controlled by boundary springs and dampers:
\begin{eqnarray}\label{1}
\left\{ \begin{array}{ll}
\rho(x)u_{tt}+\mu(x)u_{t}+\left(r(x)u_{xx}\right)_{xx}=0,\, (x,t) \in \Omega_T, \\ [4pt]
    u(x,0)=u_0(x), ~ u_t(x,0)=u_1(x),\, x\in (0,\ell),\\ [4pt]
	u(0,t)=u_{x}(0,t)=0,~ \left(-r(x)u_{xx}\right)_{x=\ell}=k_r u_x(\ell,t)+k_a u_{xt}(\ell,t), \\ [4pt]
\qquad \qquad \left(-\left(r(x)u_{xx}\right)_{x}\right)_{x=\ell}=-k_d\, u(\ell,t)-k_v u_t(\ell,t),~t\in [0,T],
\end{array} \right.
\end{eqnarray}
where $\Omega_T=(0,\ell)\times(0,T)$, $\ell>0$ is the length of the beam and $T>0$ is the final time.

Here and below, $u(x,t)$ is the vertical displacement, $r(x):=E(x)I(x)>0$ is the flexural rigidity (or bending stiffness) of the beam while $E(x)>0$ is the elasticity modulus and $I(x)>0$ is the moment of inertia of the cross section. The non-negative coefficient $\mu(x)$ represents the viscous external damping. Furthermore, the following variables have engineering meanings: $u_t(x,t)$, $u_x(x,t)$, $u_{xt}(x,t)$, $u_{xx}(x,t)$, $-\left(r(x)u_{xx}\right)$ and $-\left(r(x)u_{xx}\right)_{x}$ are the velocity, rotation, angular velocity, curvature, moment and shear force, respectively \cite{Inman:2015}. The nonnegative constants $k_r,k_d \ge 0$ and $k_a,k_v \ge 0$ represent the boundary springs and dampers, respectively.

The first boundary condition $\left(-r(x)u_{xx}\right)_{x=\ell}=k_r u_x(\ell,t)+k_v u_{xt}(\ell,t)$ at $x=\ell$ means the control resulting from the linear combination of rotation and angular velocity, and the second boundary condition $\left(-\left(r(x)u_{xx}\right)_{x}\right)_{x=\ell}=-k_d\, u(\ell,t)-k_v u_t(\ell,t)$ means the control resulting from the linear combination of displacement and  velocity. In this context, the above constants $k_r,k_d,k_a, k_v$ are defined also as the boundary controls. It should be emphasized in almost all flexible structures modeled by the Euler-Bernoulli equation, one or another special case of these boundary conditions is used (see \cite{Chen:1987a, Chen:1987b, Guo:2001, Guo:2002, Guo-Wang:2005, Krall:1989, Lazzari:2012, Toure:2015, Wang:2005} and references therein). Namely, it is shown in \cite{Guo:2001} that the generalized eigenvalues of the simplest undamped Euler-Bernoulli equation $u_{tt}-u_{xxxx}=0$ with boundary linear feedback control $u_{xx}(\ell,t)=-k_a u_{xt}(\ell,t)$, $u_{xxx}(\ell,t)=k_v u_{t}(\ell,t)$, form a Riesz basis in the state Hilbert space, which leads to exponential stability. Furthermore, for the case when $k_r=k_a= k_d=0$ and $\mu(x)=0$, the Riesz basis property and the stability of the system was studied in \cite{Guo:2002}. The same issues were studied in \cite{Toure:2015} for the system (\ref{1}) with $k_r=k_a=0$. Other simplified versions of the model governed by (\ref{1}) have been used for the mast control system in the Control of Flexible Structures Program of NASA \cite{Chen:1987a, Chen:1987b, Krall:1989}.  In \cite{Chen:1987a}, the authors examine and prove for the first time that there is exponential stability in the situation where only rotational damping is present at the extreme of a cantilever beam, with applications to long flexible structures that are modeled by the Euler-Bernoulli equation. In \cite{Chen:1987b}, the often encountered configuration in engineering practice, in which there is a finite number of serially connected beams, is analysed. In it, it is, the problem of proving uniform exponential stability when one  damper is positioned at the extremes of the composite structure, or at some intermediate interconnecting node. This problem is of great interest for structural engineers.

In all the above cited works, semigroup approach was used to obtain the Riesz basis property of the eigenfunctions, which is one of the fundamental properties of a linear vibrating system. It is well known that for such a Riesz system, the stability is usually determined by the spectrum of the associated operator. However, in the exponential stability estimate $\mathcal{E}(t) \le M e^{-\omega t} \mathcal{E}(0)$, obtained in the above mentioned studies, the relationship of the decay rate parameter $\omega >0$ with the physical and geometric parameters of the beam, including the damping coefficient and the boundary dampers, has not been determined. In addition, it does not seem possible to obtain this relationship anyway, due to the methods used in these studies.

\section{Energy identitity and dissipativity of system (\ref{1})}

We assume that the inputs in (\ref{1}) satisfy the following basic conditions:
\begin{eqnarray} \label{2}
\left \{ \begin{array}{ll}
\rho,\mu, \in L^{\infty}(0,\ell),~r \in H^2(0,\ell),\\ [3pt]
u_0\in H^2(0,\ell),\,u_1\in L^2(0,\ell),\\ [3pt]
0<\rho_0\leq\rho(x)\leq \rho_1,~ 0< r_0 \leq r(x) \leq r_1,\\ [3pt]
0\leq \mu_0\leq \mu(x)\leq \mu_1,~x \in (0,\ell)\\  [3pt]
k_r,k_a,k_d,k_v\ge 0, ~~k_a+k_v+\mu_0>0.
\end{array} \right.
\end{eqnarray}

Following the procedure described in \cite{Hasanov-Romanov:2021, Sakthivel-Hasanov:2023}, one can prove that under conditions (\ref{2}), there exists a regular weak solution $u\in L^2(0,T; H^4(0,\ell))$, $u_t\in L^2(0,T; \mathcal{V}^2(0,\ell))$  with $u_{tt}\in L^2(0,T;L^2(0,\ell))$ of problem (\ref{1}), where $\mathcal{V}^2(0,\ell):=\{v\in H^2(0,\ell):\, v(0)=v'(0)=0 \}$.

\begin{proposition}\label{Proposition-1}
Assume that conditions (\ref{2}) are satisfied. Then the following energy identity holds:
\begin{eqnarray}\label{3}
\mathcal{E}(t) + \int_0^t \int_0^\ell \mu(x) u_{\tau}^2 (x,\tau) dx d \tau  \qquad \qquad    \qquad \qquad \qquad \qquad \nonumber \\ [1pt]
\qquad  =\mathcal{E}(0)-k_a \int_0^t u_{x\tau}^2(\ell,\tau) d \tau -k_v \int_0^t u_{\tau}^2(\ell,\tau) d \tau,\,t\in[0,T],
\end{eqnarray}
where
\begin{eqnarray}\label{4}
\mathcal{E}(t)=\frac{1}{2} \int_0^\ell \left [ \rho(x) u^2_{t}(x,t)+r(x) u^2_{xx}(x,t)\right ] dx \qquad \quad \nonumber \\ [1pt]
+\frac{1}{2}\,k_r u_{x}^2(\ell,t) +\frac{1}{2}\,k_d \, u^2 (\ell,t),~t\in[0,T],
\end{eqnarray}
is the total energy of system (\ref{1}) and
\begin{eqnarray}\label{5}
\mathcal{E}(0)=\frac{1}{2} \int_0^\ell \left [ \rho(x)\left ( u_{1}(x)\right)^2 +
r(x) \left ( u''_{0}(x)\right)^2\right ] dx \qquad \quad \nonumber \\ [1pt]
+\frac{1}{2}\, k_r \left ( u'_{0}(\ell)\right)^2+\frac{1}{2}\, k_d \left ( u_{0}(\ell)\right)^2
\end{eqnarray}
is the initial value of the total energy.
\end{proposition}
\textbf{Proof.} Multiplying both sides of equation (\ref{1}) by $u_t(x,t)$, integrating it over $\Omega_t:=(0,\ell)\times (0,t)$, using the identity
\begin{eqnarray}\label{6}
(r(x)u_{xx})_{xx} u_t =  [(r(x)u_{xx})_x u_t-r(x)u_{xx} u_{xt}]_x+ \frac{1}{2}\left (r(x)u_{xx}^2\right )_t,
\end{eqnarray}
we obtain the following integral identity:
\begin{eqnarray*}
\frac{1}{2} \int_0^t\int_0^\ell \left (\rho(x) u_{\tau}^2\right )_{\tau}dx\,d\tau +\frac{1}{2} \int_0^t \int_0^\ell \left (r(x)u_{xx}^2\right )_{\tau}dx\,d\tau \qquad \qquad  \\ [1pt]
\qquad +\int_0^t \left ((r(x)u_{xx})_x u_{\tau}-r(x)u_{xx} u_{x\tau} \right)_{x=0}^{x=\ell} d\tau +\int_0^t \int_0^\ell \mu(x) u_{\tau}^2 dx d \tau=0,
\end{eqnarray*}
for all $t \in (0,T]$. Using here the initial and boundary conditions (\ref{2}), we obtain:
 \begin{eqnarray*}
\frac{1}{2} \int_0^\ell \left [\rho(x) u^2_{t}+ r(x) u_{xx}\right ]dx
+\frac{1}{2}\,k_r u_{x}^2(\ell,t) +\frac{1}{2}\,k_d \, u^2 (\ell,t)\qquad \qquad \\ [1pt]
+\int_0^t \int_0^\ell \mu(x) u_{\tau}^2 dx d \tau+k_a \int_0^t u_{x\tau}^2(\ell,\tau) d \tau +k_v \int_0^t u_{\tau}^2(\ell,\tau) d \tau \qquad  \\ [1pt]
-\frac{1}{2} \int_0^\ell \left [\rho(x) \left (u_1(x)\right )^2+ r(x) \left (u''_0(x)\right )^2\right ]dx
-\frac{1}{2} k_r \left ( u'_{0}(\ell)\right)^2-\frac{1}{2} k_d \left ( u_{0}(\ell)\right)^2=0,
\end{eqnarray*}
for all $t \in (0,T]$. This leads to (\ref{3}) with (\ref{4}) and (\ref{5}).\hfill$\Box$

\medskip

Identity (\ref{3}) show that the increase in the damping  parameters $k_a,\,k_v \ge 0$ causes the energy $\mathcal{E}(t)$ to decrease, Furthermore, from formula (\ref{4}) it follows that the increase of the spring parameters $k_r,\,k_d \ge 0$ causes the energy to increase.

\begin{proposition}\label{Proposition-2}
If conditions (\ref{2}) are met, the formula below gives the rate at which the total energy decreases.
\begin{eqnarray}\label{7}
\frac{d \mathcal{E}(t)}{dt} =-\int_0^\ell \mu(x) u^2_{t}dx- k_a u_{xt}^2 (\ell,t) - k_v u_{t}^2(\ell,t),\, t\in (0,T).
\end{eqnarray}
\end{proposition}
\textbf{Proof.} In view of formula (\ref{4}) we have:
\begin{eqnarray*}
\frac{d \mathcal{E}(t)}{dt}= \int_0^\ell \left [ \rho(x) u_{t}u_{tt}+r(x) u_{xx}u_{xxt}\right ] dx \qquad \qquad \qquad \qquad\\ [1pt]
\qquad \qquad  \qquad +k_r u_{x} (\ell,t) u_{xt} (\ell,t) + k_d u(\ell,t) u_{t}(\ell,t),~t\in[0,T].
\end{eqnarray*}
Use here the (formal) identity $\rho(x)u_{tt}=-\mu(x)u_{t}-\left(r(x)u_{xx}\right)_{xx}$ to get
\begin{eqnarray}\label{8}
\frac{d \mathcal{E}(t)}{dt}=- \int_0^\ell \mu(x) u^2_{t}dx- \int_0^\ell \left(r(x)u_{xx}\right)_{xx} u_{t} \,dx+\int_0^\ell r(x) u_{xx}u_{xxt}\,dx \nonumber \\ [1pt]
\qquad \qquad  \qquad +k_r u_{x} (\ell,t) u_{xt} (\ell,t) + k_d u(\ell,t) u_{t}(\ell,t),~t\in[0,T].
\end{eqnarray}
In the second right hand side integral, we employ the identity
\begin{eqnarray*}
-\int_0^\ell \left(r(x)u_{xx}\right)_{xx} u_{t} \,dx =-\int_0^\ell r(x) u_{xx}u_{xxt}dx-k_d\, u(\ell,t) u_{t} (\ell,t)
 \nonumber\\ [1pt]
\qquad \qquad  - k_v u^2_{t}(\ell,t) -k_r u_x(\ell,t) u_{xt} (\ell,t) - k_a u^2_{xt}(\ell,t),\, t\in[0,T],
\end{eqnarray*}
which holds due to the boundary conditions in (\ref{1}). Substituting this identity in (\ref{8}) we arrive at the required formula (\ref{7}). \hfill$\Box$

\begin{corollary}\label{Corollary-1}
Integrating (\ref{7}) over $(0,t)$, $t \in (0,T]$ we obtain the same energy identity (\ref{3}) rewritten in the following form:
\begin{eqnarray}\label{9}
\mathcal{E}(0)-\mathcal{E}(t)=\mathfrak{j}_\mu (t)+\mathfrak{j}_a (t) +\mathfrak{j}_v (t),~t\in[0,T],
\end{eqnarray}
where
\begin{eqnarray} \label{10}
\left. \begin{array}{ll}
\displaystyle \mathfrak{j}_\mu (t): =  \int_0^t\int_0^\ell \mu(x) u^2_{\tau}(x,\tau) dx d\tau,\\ [12pt]
\displaystyle \mathfrak{j}_a(t): = k_a \int_0^t  u_{x\tau}^2 (\ell,\tau)d\tau,~~~~~\mathfrak{j}_v (t): = k_v \int_0^t u_{\tau}^2(\ell,\tau)d\tau,~t\in[0,T].
\end{array} \right.
\end{eqnarray}

In particular,
\begin{eqnarray*}
\mathcal{E}(t) \le \mathcal{E}(0),~t\in[0,T],
\end{eqnarray*}
that is, the energy of the system (\ref{1}) is dissipating with time.
\end{corollary}

The above formula (\ref{7}) is a clear expression of the effect of the damping parameters $\mu(x)$, $k_a$ and $k_v$ on the rate of decrease of the total energy. In addition, the energy identity (\ref{9}) shows the degree of influence of these damping factors on the difference between the initial value $\mathcal{E}(0)$ of the total energy and the value $\mathcal{E}(t)$ of this energy at the time instant $t\in (0,T]$, through the integrals $\mathfrak{j}_\mu (t)$, $\mathfrak{j}_a (t)$ and $\mathfrak{j}_v(t)$ defined in (\ref{10}).

\section{Energy decay estimate for system (\ref{1})}

Introduce the auxiliary function:
\begin{eqnarray}\label{11}
\mathcal{J}(t)= \int_0^\ell \rho(x) u u_{t}dx+\frac{1}{2} \int_0^\ell \mu(x) u^2dx \qquad \qquad \qquad \nonumber\\ [2pt]
\qquad +\frac{1}{2} \,k_a u_{x}^2 (\ell,t) +\frac{1}{2}\, k_v u^2(\ell,t),\, t\in[0,T],
\end{eqnarray}
containing all damping parameters.

We prove the formula
\begin{eqnarray}\label{12}
\frac{d \mathcal{J}(t)}{dt}= 2 \int_0^\ell \rho(x) u_t^2dx -2\mathcal{E}(t),~t\in[0,T],
\end{eqnarray}
which shows the relationship between the auxiliary function $\mathcal{J}(t)$ and the energy function $\mathcal{E}(t)$ introduced in (\ref{4}).

Taking the derivative of the function $\mathcal{J}(t)$ with respect to the time variable and using then the (formal) identity $\rho(x)u_{tt}+\mu(x)u_{t}=-\left(r(x)u_{xx}\right)_{xx}$ as above, we obtain:
\begin{eqnarray*}
\frac{d \mathcal{J}(t)}{dt}= \int_0^\ell \rho(x) u^2_{t}dx
- \int_0^\ell \left(r(x)u_{xx}\right)_{xx}\,u \,dx \qquad \qquad \qquad  \\ [2pt]
 +k_a u_{x} (\ell,t) u_{xt} (\ell,t) +k_v u(\ell,t) u_{t} (\ell,t),~t\in[0,T].
\end{eqnarray*}
We employ here the following identity:
\begin{eqnarray*}
-\int_0^\ell \left(r(x)u_{xx}\right)_{xx} u \,dx =
-\int_0^\ell r(x) u^2_{xx}dx \qquad \qquad \qquad \qquad \qquad \qquad \nonumber\\
\qquad  -k_d u^2(\ell,t)-k_v u(\ell,t) u_{t}(\ell,t)
-k_r u^2_{x}(\ell,t)-k_a u_{x}(\ell,t) u_{xt} (\ell,t),
\end{eqnarray*}
$t\in[0,T]$. This yields:
\begin{eqnarray*}
\frac{d \mathcal{J}(t)}{dt}= \int_0^\ell \rho(x) u^2_{t}dx
- \int_0^\ell r(x) u^2_{xx} dx \qquad \qquad \qquad  \\ [2pt]
-k_d\, u^2(\ell,t)-k_r u^2_{x}(\ell,t),~t\in[0,T].
\end{eqnarray*}
With definition (\ref{4}) this implies the desired formula (\ref{12}). \hfill$\Box$

\begin{proposition}\label{Proposition-3}
Under conditions (\ref{2}), the energy function $\mathcal{E}(t)$ introduced in (\ref{4}) serves as lower and upper bounds to the auxiliary function $\mathcal{J}(t)$ introduced in (\ref{11}), that is
\begin{eqnarray}\label{13}
-\beta_0 \,\mathcal{E}(t) \le \mathcal{J}(t) \le \beta_1\, \mathcal{E}(t), ~t\in[0,T],
\end{eqnarray}
where
\begin{eqnarray}\label{14}
\left. \begin{array}{ll}
\displaystyle \beta_1=\beta_0 \left [1+  \frac{1}{\sqrt{\rho_1 r_0}}\left (\frac{\ell^2}{2}\,\mu_1+\frac{2}{\ell}\,k_a
+ \ell\,k_v \right )\right ],\\ [14pt]
\displaystyle \beta_0 =\frac{\ell^2}{2}\,\sqrt{\frac{\rho_1}{r_0}}\,. \\
\end{array} \right.
\end{eqnarray}
\end{proposition}
\textbf{Proof.} First we estimate the first right hand side integral in (\ref{11}).  To this end, we employ the $\varepsilon$-inequality
\begin{eqnarray*}
\left \vert \int_0^\ell \rho(x) u u_{t}dx\right \vert \le \frac{\varepsilon}{2}\, \int_0^\ell \rho(x) u_t^2dx + \frac{1}{2\varepsilon}\, \int_0^\ell \rho(x) u^2dx,
\end{eqnarray*}
with the inequality
\begin{eqnarray*}
\int_0^\ell \rho(x) u^2dx \le \frac{\ell^4 \rho_1}{4 r_0} \int_0^\ell r(x)u_{xx}^2dx,
\end{eqnarray*}
to estimate the second right-hand side integral in above inequality. Choosing then the parameter $\varepsilon>0$ from the condition $\varepsilon/2=\ell^4 \rho_1/(8 r_0\,\varepsilon)$ as
\begin{eqnarray*}
\varepsilon= \frac{\ell^2}{2}\,\sqrt{\frac{\rho_1}{r_0}}\,,
\end{eqnarray*}
we obtain the following estimate:
\begin{eqnarray}\label{15}
\left \vert \int_0^\ell \rho(x) u u_{t}dx\right \vert \le
\frac{\ell^2}{4}\,\sqrt{\frac{\rho_1}{r_0}}\left [ \int_0^\ell \rho(x) u_t^2dx + \int_0^\ell r(x) u_{xx}^2dx \right ].
\end{eqnarray}

For other right-hand side terms in formula (\ref{11}) for the auxiliary function  $\mathcal{J}(t)$ we use the following inequalities:
\begin{eqnarray}\label{16}
\left. \begin{array}{ll}
\displaystyle \frac{1}{2} \int_0^\ell \mu(x) u^2 (x,t)dx \le
\frac{\ell^4 \mu_1}{16 r_0} \int_0^\ell r(x)u_{xx}^2dx,\\ [10pt]
\displaystyle \frac{1}{2}\,k_a u^2_{x}(\ell,t)\le
\frac{\ell\, k_a}{2 r_0} \int_0^\ell r(x)u_{xx}^2dx,\\ [10pt]
\displaystyle \frac{1}{2}\,k_v u^2(\ell,t)\le \frac{\ell^3 k_v}{4 r_0} \int_0^\ell r(x)u_{xx}^2dx,~t\in (0,T).
\end{array} \right.
\end{eqnarray}

Taking into account (\ref{15}) and (\ref{16}) in (\ref{11}) we arrive at the following estimate
\begin{eqnarray*}
&& \mathcal{J}(t) \le \frac{1}{2}\, \beta_0 \left \{ \int_0^\ell \rho (x)u^2_{t}dx \right.  \\ [2pt]
&& \qquad \qquad +\, \left. \left [1+ \frac{\ell^2}{2\,\sqrt{\rho_1r_0}}\,\mu_1 +
\frac{2}{\ell\,\sqrt{\rho_1r_0}}\,k_a +\frac{\ell}{\sqrt{\rho_1r_0}}\,k_v
\right ] \int_0^\ell r(x)u^2_{xx}dx  \right \},
\end{eqnarray*}
 which leads to the upper bound
\begin{eqnarray}\label{17}
\mathcal{J}(t) \le \beta_1 \,\mathcal{E}(t),\, t\in[0,T], ~\beta_1 >0,
\end{eqnarray}
with $\beta_0, \beta_1>0$ introduced in (\ref{13}).

To find the lower bound for the auxiliary function $\mathcal{J}(t)$, we use again inequality (\ref{15}) in (\ref{11}) to conclude that

\begin{eqnarray*}
\mathcal{J}(t) \ge - \frac{1}{2}\, \beta_0 \left \{ \int_0^\ell \rho (x)u^2_{t}dx +\int_0^\ell r(x)u^2_{xx}dx \right \}\ + \qquad \quad \qquad \quad \qquad \qquad \\ [2pt]
 \qquad \quad \qquad \quad  \frac{1}{2} \int_0^\ell \mu(x) u^2 (x,t)dx\ +\  \frac{1}{2} \,k_a u_{x}^2 (\ell,t) +\frac{1}{2}\, k_v u^2(\ell,t), ~t\in[0,T],
\end{eqnarray*}
This leads to
\begin{eqnarray}\label{18}
\mathcal{J}(t) \ge -\beta_0\, \mathcal{E}(t),~t\in [0,T].
\end{eqnarray}

Thus, (\ref{17}) and (\ref{18}) imply the required lower and upper bounds (\ref{14}). \hfill$\Box$

\begin{remark}\label{Remark-1}
The constants $\beta_0,\beta_1>0$ depend only on the geometric and physical parameters of a beam introduced in (\ref{2}), as formulas (\ref{14}) show.
\end{remark}

\medskip

To establish the uniform energy decay estimate, we introduce the Lyapunov function:
\begin{eqnarray}\label{19}
\mathcal{L}(t)=\mathcal{E}(t)+\lambda\mathcal{J}(t),\, t\in[0,T],
\end{eqnarray}
where $\mathcal{E}(t)$ and $\mathcal{J}(t)$ are the energy function and the auxiliary function  introduced in (\ref{4}) and (\ref{11}), respectively, and $\lambda>0$ is the penalty term.

\begin{theorem}\label{Theorem-1}
Assume that conditions (\ref{2}) are satisfied. Then system (\ref{1}) is exponentially stable for any nonnegative values of the boundary spring and damper constants $k_r, k_d, k_a, k_v\ge 0$. That is, there are the constants
\begin{eqnarray}\label{20}
\displaystyle M_d= \frac{1+  \beta_1 \lambda}{1-  \beta_0 \lambda }~, ~
 \sigma=\frac{2 \lambda}{1+\beta_1 \lambda}~,
\end{eqnarray}
with
\begin{eqnarray}\label{21}
0<\lambda <\min (1/ \beta_0,\, \mu_0/(2\rho_1)).
\end{eqnarray}
such that the energy $\mathcal{E}(t)$ of system (\ref{1}) satisfies the following estimate:
\begin{eqnarray}\label{22}
\mathcal{E}(t)\le M_d\,  e^{-\sigma t}\, \mathcal{E}(0),~t\in[0,T],
\end{eqnarray}
where $\mu_0,\rho_1>0$ and $\beta_0>0$ are the constants introduced in (\ref{2}) and (\ref{14}), respectively, and $\mathcal{E}(0)>0$ is the initial energy defined in (\ref{5}).
\end{theorem}
\textbf{Proof.} In view of (\ref{14}) we have:
\begin{eqnarray}\label{23}
\left (1- \lambda \, \beta_0 \right ) \mathcal{E}(t) \le \mathcal{L}(t) \le  \left (1+ \lambda \, \beta_1 \right )\mathcal{E}(t),\, t\in[0,T].
\end{eqnarray}
In such a circumstance, we assume that the penalty term satisfies the following conditions:
\begin{eqnarray}\label{24}
0<\lambda <1/ \beta_0,~\beta_0>0.
\end{eqnarray}

Differentiating $\mathcal{L}(t)$ with respect to the variable $t\in (0,T)$ and taking formulas (\ref{7}) and (\ref{12}) into account, we obtain:
\begin{eqnarray}\label{25}
\frac{d \mathcal{L}(t)}{dt}+2 \lambda \mathcal{E}(t)= -\int_0^\ell \left [\mu(x)-2\lambda \rho(x)\right ] u_t^2dx
 \qquad \qquad \quad \nonumber\\ [2pt]
\qquad \qquad \qquad -k_a u^2_{xt}(\ell,t)-k_v u^2_{t}(\ell,t),~t\in[0,T].
\end{eqnarray}
We require that $\mu(x)-2\lambda \rho(x)>0$. Since $\mu(x) -2\lambda \rho(x) \ge \mu_0-2\lambda \rho_1$, the sufficient condition for this is the condition
\begin{eqnarray}\label{26}
\lambda< \mu_0/(2\rho_1).
\end{eqnarray}
With (\ref{24}) this implies that the penalty term should satisfy conditions (\ref{21}).
Then from (\ref{25}) we deduce that
\begin{eqnarray}\label{27}
\frac{d \mathcal{L}(t)}{dt}+2 \lambda \mathcal{E}(t)<0,~t\in[0,T].
\end{eqnarray}
With the inequality $ \mathcal{E}(t) \ge \mathcal{L}(t)/ \left (1+ \lambda \, \beta_1 \right )$ this yields:
\begin{eqnarray*}
\frac{d \mathcal{L}(t)}{dt}+\frac{2 \lambda}{1+\beta_1 \lambda} \, \mathcal{L}(t)<0,~t\in[0,T].
\end{eqnarray*}
Solving this inequality we find:
\begin{eqnarray*}
\mathcal{L}(t)\le  e^{-\sigma t} \mathcal{L}(0),~t\in[0,T].
\end{eqnarray*}
This yields the required estimate (\ref{22}) with the constants $M_d,\sigma >0$ introduced in (\ref{20}). \hfill$\Box$

\begin{remark}\label{Remark-2}
In view of formulas (\ref{14}), the decay rate parameter $\sigma>0$ in the energy estimate (\ref{22}), obtained for the system governed by (\ref{1}) and controlled by boundary springs and dampers, clearly show the degree of influence of each of the damping parameters $\mu(x)$, $k_a, k_v \ge 0$ in the dissipative boundary conditions on the energy decay.
\end{remark}

\section{Some special cases}

Special cases of the general system (\ref{1}) described above are very common in practical applications of structures containing beam elements. In this section we deal with systems corresponding to special cases of the general system (\ref{14}) to investigate the influence of each damping factor.

\subsection{A cantilever beam fixed at one end and free at other}
\medskip

Consider the simplest case when $k_r=k_a=k_d=k_v=0$ of system (\ref{1}), i.e. without the dissipative boundary conditions:
\begin{eqnarray}\label{1a}
\left\{ \begin{array}{ll}
\rho(x)u_{tt}+\mu(x)u_{t}+\left(r(x)u_{xx}\right)_{xx}=0,\, (x,t) \in \Omega_T, \\ [4pt]
    u(x,0)=u_0(x), ~ u_t(x,0)=u_1(x),\, x\in (0,\ell),\\ [4pt]
	u(0,t)=u_{x}(0,t)=0,~ \left(-r(x)u_{xx}\right)_{x=\ell}=0, \\ [4pt]
\qquad \qquad \qquad \qquad \qquad \left(-\left(r(x)u_{xx}\right)_{x}\right)_{x=\ell}=0,~t\in [0,T].
\end{array} \right.
\end{eqnarray}
This is an initial boundary value problem for the damped cantilever beam.

The exponential stability result for system (\ref{1a}) directly follows from the results given in (\ref{20})-(\ref{22}),
\begin{eqnarray}\label{20a}
\left\{ \begin{array}{ll}
\mathcal{E}(t)\le M_0\,  e^{-\sigma_0 t}\, \mathcal{E}(0),~t\in[0,T], \\ [6pt]
\displaystyle M_0= \frac{1+  \beta_1 \lambda}{1-  \beta_0 \lambda }\,, ~
 \sigma_0=\frac{2 \lambda}{1+\beta_1 \lambda}~,\\ [10pt]
 \displaystyle \beta_1=\beta_0 \left [1+  \frac{\ell^2}{4\sqrt{\rho_1 r_0}}\,\mu_1\right ],~
\displaystyle \beta_0 =\frac{\ell^2}{2}\,\sqrt{\frac{\rho_1}{r_0}}\, , \\ [14pt]
 0<\lambda <\min (1/ \beta_0,\, \mu_0/(2\rho_1)),
\end{array} \right.
\end{eqnarray}
assuming $k_a=k_v=0$ in (\ref{14}), and also $k_r=k_d=0$ in (\ref{4}). That is, the energy function corresponding to system (\ref{1a}) is
 \begin{eqnarray}\label{4a}
\mathcal{E}(t)=\frac{1}{2} \int_0^\ell \left [ \rho(x) u^2_{t}(x,t)+r(x) u^2_{xx}(x,t)\right ] dx,~t\in[0,T].
\end{eqnarray}

Formulas (\ref{20a}) clearly show the nature of the influence of the viscous external damping coefficient $\mu(x)$, as a unique damping factor on the energy decay rate.

\subsection{A cantilever beam fixed at one end and attached to a spring at other}
\medskip

This case corresponds to the zero values $k_a=k_v=0$ of the boundary damping parameters, and hence to the linear spring conditions at $x=\ell$:
\begin{eqnarray}\label{1a1}
\left\{ \begin{array}{ll}
\rho(x)u_{tt}+\mu(x)u_{t}+\left(r(x)u_{xx}\right)_{xx}=0,\, (x,t) \in \Omega_T, \\ [4pt]
    u(x,0)=u_0(x), ~ u_t(x,0)=u_1(x),\, x\in (0,\ell),\\ [4pt]
	u(0,t)=u_{x}(0,t)=0,~ \left(-r(x)u_{xx}\right)_{x=\ell}=k_r u_{x}(\ell,t), \\ [4pt]
\qquad \qquad \qquad \qquad \qquad \left(-\left(r(x)u_{xx}\right)_{x}\right)_{x=\ell}=k_d\, u(\ell,t),~t\in [0,T].
\end{array} \right.
\end{eqnarray}
As in the previous case, the dissipativity of system (\ref{1a1}) is provided only by the viscous external damping given by the coefficient $\mu(x)>0$.

The same exponential stability result given in (\ref{20a}) holds for system (\ref{1a1}).
Furthermore, the energy function $\mathcal{E}(t)$ corresponding to system (\ref{1a1}) is given by the same formula (\ref{4}) which, different from formula (\ref{4a}), contains also the spring constants $k_r,k_d \ge 0$.

 \subsection{A cantilever beam fixed at one end and subjected two dampers at other}
 \medskip

 Consider the case where both spring parameters in (\ref{1}) are zero, $k_r=k_d=0$:
 \begin{eqnarray}\label{1b}
\left\{ \begin{array}{ll}
\rho(x)u_{tt}+\mu(x)u_{t}+\left(r(x)u_{xx}\right)_{xx}=0,\, (x,t) \in \Omega_T, \\ [4pt]
    u(x,0)=u_0(x), ~ u_t(x,0)=u_1(x),\, x\in (0,\ell),\\ [4pt]
	u(0,t)=u_{x}(0,t)=0,~ \left(-r(x)u_{xx}\right)_{x=\ell}=k_a u_{xt}(\ell,t), \\ [4pt]
\qquad \qquad \qquad \left(-\left(r(x)u_{xx}\right)_{x}\right)_{x=\ell}=-k_v\, u_t(\ell,t),~t\in [0,T],
\end{array} \right.
\end{eqnarray}

This is a mathematical model for the mast control system. The simplest version
 \begin{eqnarray}\label{1bb}
\left\{ \begin{array}{ll}
m\,u_{tt}+EI\,u_{xxxx}=0,\, (x,t) \in \Omega_T, \\ [4pt]
    u(x,0)=u_0(x), ~ u_t(x,0)=u_1(x),\, x\in (0,\ell),\\ [4pt]
	u(0,t)=u_{x}(0,t)=0,~ \left(-EI\,u_{xx}\right)_{x=\ell}=k_a u_{xt}(\ell,t), \\ [4pt]
\qquad \qquad \qquad \left(-EI\,u_{xxx}\right)_{x=\ell}=-k_v\, u_t(\ell,t),~t\in [0,T],
\end{array} \right.
\end{eqnarray}
of this model for the undamped Euler-Bernoulli equation with constant coefficients was first studied in \cite{Chen:1987a} within NASA's Program of Control of Flexible Structures, and then developed in  \cite{Chen:1987b}.

In this model, the meaning of the boundary conditions at $x=\ell$ is that the shear force $-EI\, u_{xxx}$ is proportional to velocity $u_t$, and the bending moment $-EI\, u_{xx}$ is negatively proportional to angular velocity $u_{xt}$, while the values of the boundary dampers $k_a, k_v \ge 0$ play the role of the proportionality factors. Thus, the rate feedback laws at $x=\ell$ reflect basic features of mast control systems with bending and torsion rate control.

The uniform exponential stability result $\mathcal{E}(t)\le K\,  e^{-\mu\, t}\, \mathcal{E}(0)$ for the energy of vibration of the beam governed by system (\ref{1bb}) was proved in \cite{Chen:1987b}. However, the constants $K, \mu>0$ are not related to either physical or boundary damping parameters. Therefore, from this estimate, it is impossible to reveal the degree of influence of these parameters on energy decay.

The results given in (\ref{20})-(\ref{22}), with the same constants $\beta_0,\beta_1>0$  introduced in (\ref{14}),  are valid also for system (\ref{1b}). However, in the case of $\mu(x)=0$, the sufficient condition (\ref{26}) for ensuring the inequality (\ref{27}) cannot be given over the coefficient $\mu(x)$. As a consequence, the above results can not be used for system (\ref{1bb}) with undamped Euler-Bernoulli equation.

\begin{theorem}\label{Theorem-2}
Assume that conditions (\ref{2}) are satisfied and $\mu(x)=0$. Suppose, in addition that
 \begin{eqnarray}\label{ad-1b}
\left. \begin{array}{ll}
u^2_{xt}(\ell,t)+u^2_{t}(\ell,t)>0,~\mbox{for\,all}~t\in [0,T].
\end{array} \right.
\end{eqnarray}
Then system (\ref{1bb}) is exponentially stable:
\begin{eqnarray}\label{22b}
\mathcal{E}(t)\le M_d\,  e^{-\sigma t}\, \mathcal{E}(0),~t\in[0,T],
\end{eqnarray}
where
\begin{eqnarray}\label{20bb}
\left. \begin{array}{ll}
\displaystyle M_d= \frac{1+  \beta_1 \lambda}{1-  \beta_0 \lambda }~, ~
 \sigma=\frac{2 \lambda}{1+\beta_1 \lambda}~,\\ [12pt]
 \displaystyle \beta_1=\beta_0 \left [1+\frac{1}{\sqrt{m\,EI}} \left(\frac{2}{\ell}\,k_a + \ell\,k_v \right ) \right ],~ \beta_0 =\frac{\ell^2}{2}\,\sqrt{\frac{m}{EI}},\\ [14pt]
\displaystyle 0<\lambda <\min \left ( \frac{1}{\beta_0},~\frac{\inf_{[0,T]}\left [k^2_a u^2_{xt}(\ell,t)+k^2_v u^2_{t}(\ell,t)\right ]}{2m \Vert u_t \Vert_{L^{\infty}(0,T;L^2(0,\ell))}^2} \, \right).
 \end{array} \right.
\end{eqnarray}
\end{theorem}
\textbf{Proof.} This theorem is proved in the similar way as the previous theorem, one only needs to derive the similar inequality for the Lyapunov function $\mathcal{L}(t)$, through the boundary damping parameters $k_a,k_v\ge 0$. To this end, we use the following analogue
\begin{eqnarray}\label{25b}
\frac{d \mathcal{L}(t)}{dt}+2 \lambda \mathcal{E}(t)= 2\lambda\,m \int_0^\ell u_t^2dx
-k_a u^2_{xt}(\ell,t)-k_v u^2_{t}(\ell,t),\,t\in[0,T]\quad
\end{eqnarray}
of formula (\ref{25}) of the Lyapunov function, which corresponds to system  (\ref{1bb}). We require that
\begin{eqnarray*}
2\lambda\,m \int_0^\ell u_t^2dx
-k_a u^2_{xt}(\ell,t)-k_v u^2_{t}(\ell,t)<0,\,t\in[0,T]\quad
\end{eqnarray*}
Evidently, for the penalty term $\lambda>0$ satisfying the third condition of (\ref{20bb}), the above inequality holds for all $t\in[0,T]$. This implies inequality (\ref{27}).
The uniform exponential decay estimate (\ref{22b}) is obtained from this inequality in the same way as in the proof of Theorem \ref{Theorem-1}.
\hfill$\Box$

\section{Numerical results}

In many cases, especially with variable coefficients, it is not easy to obtain an analytical solution for given problem (\ref{1}).  Since the demand for finding energy function in (\ref{4}) involves $u$ and its derivatives, these quantities should be calculated by an efficient numerical technique. In the next part, we first briefly summarize a robust one which is known the Method of Lines (MOL) approach that has been used successfully in many previous studies related to Euler-Bernoulli beam equations with the classical boundary conditions \citep{AH:OB2016}-\citep{Hasanov-Kawano-Baysal:2023}. Then we both show the implementation of this method to the considered problem (\ref{1}) and  demonstrate its high accuracy performance.

\subsection{Method of Lines Approach for the Numerical Solution of (\ref{1})}

The MOL approach is based on two-stage decomposition principle for (\ref{1}); first, a semi-discrete formula is obtained from the variational formulation by Finite Element Method (FEM) with Hermite cubic shape functions, then the full discretization is generated by the second order appropriate time integrators. At the end of this process an algebraic system is obtained which is simple to solve. This technique is commonly employed, particularly in the case of dynamical multi-dimensional phenomena.

Assume the finite dimensional space $V_h\subset\mathcal{V}^2(0,\ell)$ spanned by the Hermite cubic shape functions $\{\psi_i\}_{i=1}^{2M}$ by uniformly discretizing spatial domain $0=x_1<x_2<\cdots<x_M=\ell$ (where $h=\ell/(M-1)$). Consider the following semi-discrete Galerkin approximation of the problem (\ref{1}).\\

\noindent \emph{For all $t\in(0,T]$, find $u_h(\cdot,t)\in  V_h$ such that $\forall v_h\in V_h$,
\begin{eqnarray}\label{semidisc}
\left\{\begin{array}{ll}
(u_{h,tt}(\cdot,t),v_h)+(\mu(\cdot)u_{h,t}(\cdot,t),v_h)+a(u_h(\cdot,t),v_h) =\\
\hskip .5cm -[k_d\, u_h(\ell,t)+k_v u_{h,t}(\ell,t)]v_h(\ell) -[(k_r\, u_{h,x}(\ell,t)+k_a u_{h,xt}(\ell,t)])v_{h,x}(\ell),\\
u_h(x,0)=0,~~ u_{h,t}(x,0)=0.
\end{array}\right.
\end{eqnarray}}
Here $u_h(x,t)$  is the finite element approximation of the weak solution of (\ref{1}) and the symmetric bilinear functional $a:H^2(0,\ell)\times H^2(0,\ell)\rightarrow \mathbb{R}$ is defined by $a(w,v):=(r(\cdot)w_{xx},v_{xx})$.

The above second-order system of ODE can be approximately solved by using the following second-order backward finite difference approximations of $u_{h,t}(x,t_j)$ and $u_{h,tt}(x,t_j)$ with uniform temporal discretization.
$0=t_1<t_2<\cdots<t_N=T$ (where $dt=T/(N-1)$).
 \begin{eqnarray}\label{imptimeint}
&&\hspace*{-.5cm}u_{h,t}(x,t_j)\approx \partial_{t}^-U_h^j(x):=\frac{3u_h(x,t_j)-4u_h(x,t_{j-1})+u_h(x,t_{j-2})}{2\, h_t},\nonumber \\
&&\hspace*{-.5cm}u_{h,tt}(x,t_j)\approx \partial_{tt}^-U_h^j(x):=\frac{2u_h(x,t_j)-5u_h(x,t_{j-1})+4u_h(x,t_{j-2})-u_h(x,t_{j-3})}{h_t^2}.\nonumber
\end{eqnarray}
By substituting these difference quotients for $u_h(x,t)$ in the semi-discrete analogy (\ref{semidisc}), one can get the following full-discrete algebraic problem of which solution $U_h^j(x)$ is the approximate solution of (\ref{1}) at $t=t_j$ such that $U_h^j\approx u(\cdot,t_j)$.\\

\noindent \emph{For each $j=1,2,...,N,$ find $U_h^j\in V_h$ such that $\forall v_h\in V_h$,
\begin{eqnarray}\label{fulldisc}
\left\{\begin{array}{ll}
(\partial_{tt}^- U_h^j,v_h)+(\mu(\cdot)\partial_{t}^- U_h^j,v_h)+a(U_h^j,v_h)=\\
\hskip .5cm -[k_d\, U_h^j(\ell)+k_v \partial_t^-U_{h}^j(\ell)]v_h(\ell) -[(k_r\, U_{h,x}^j(\ell)+k_a \partial_t^-U_{h,x}^j(\ell)])v_{h,x}(\ell).
\end{array}\right.
\end{eqnarray}}
In order to compare numerical and exact solution on cartesian coordinates, we define ${U_h(x)}$ as linear interpolation of set of all solutions  $\{U_h^j\in V_h\}_{j=1}^N$ in temporal dimension such that  for $j=1,\, \cdots,\, N-1$,
$$U_h(x,t)|_{[t_{j},t_{j+1}]}:=\frac{t-t_{j}}{h_t}U_h^{j+1}(x)-\frac{t-t_{j+1}}{h_t}U_h^{j}(x).$$
In the next sectionin, we will test the success of this MOL technique with a problem for which we know the exact solution and develop a simple method for approximating the desired energy function.

\subsection{Test Problem}
\begin{figure}[!htb]
	\includegraphics[width=6.8cm,height=6.8cm]{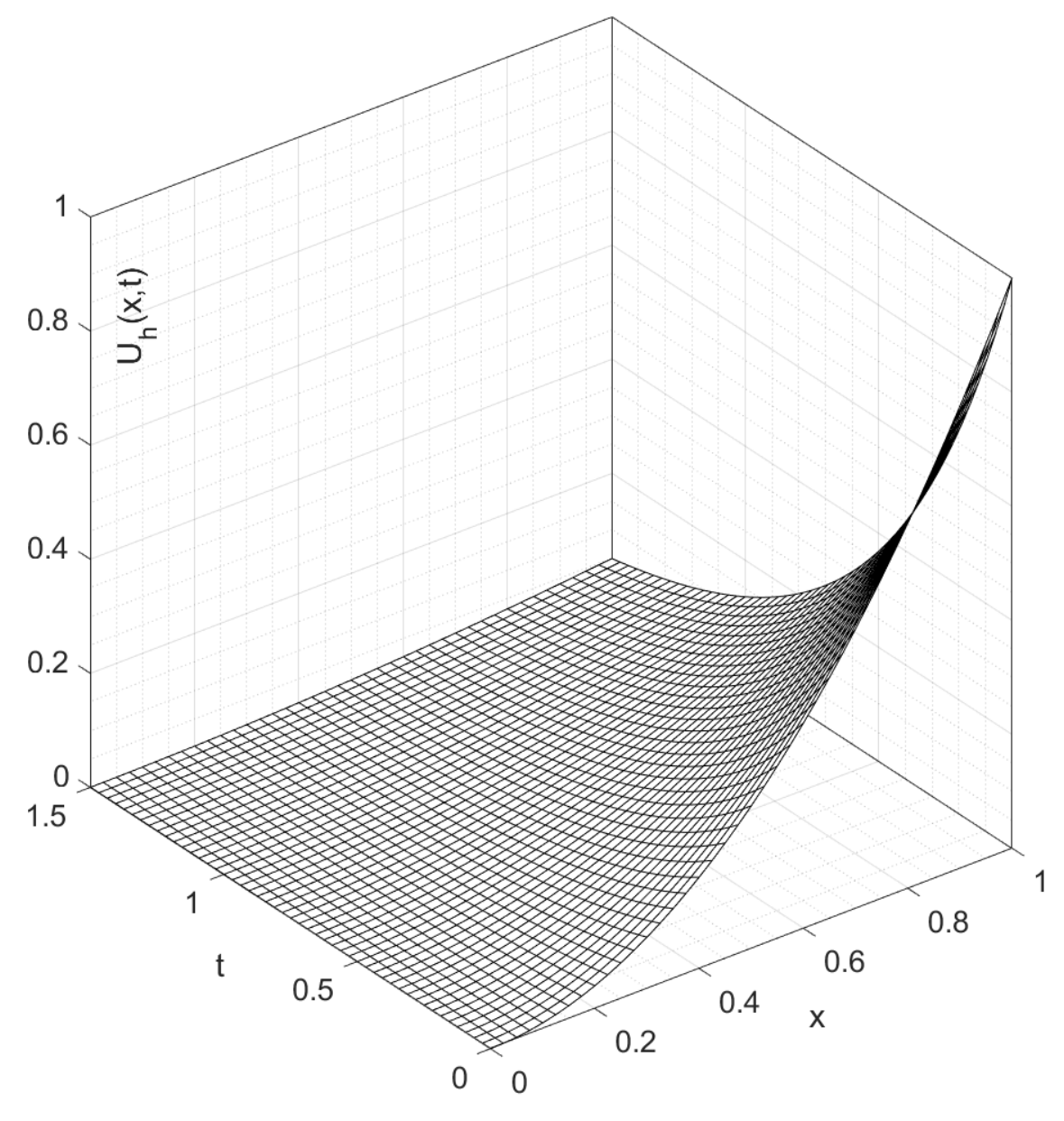}~~\includegraphics[width=6.8cm,height=6.8cm]{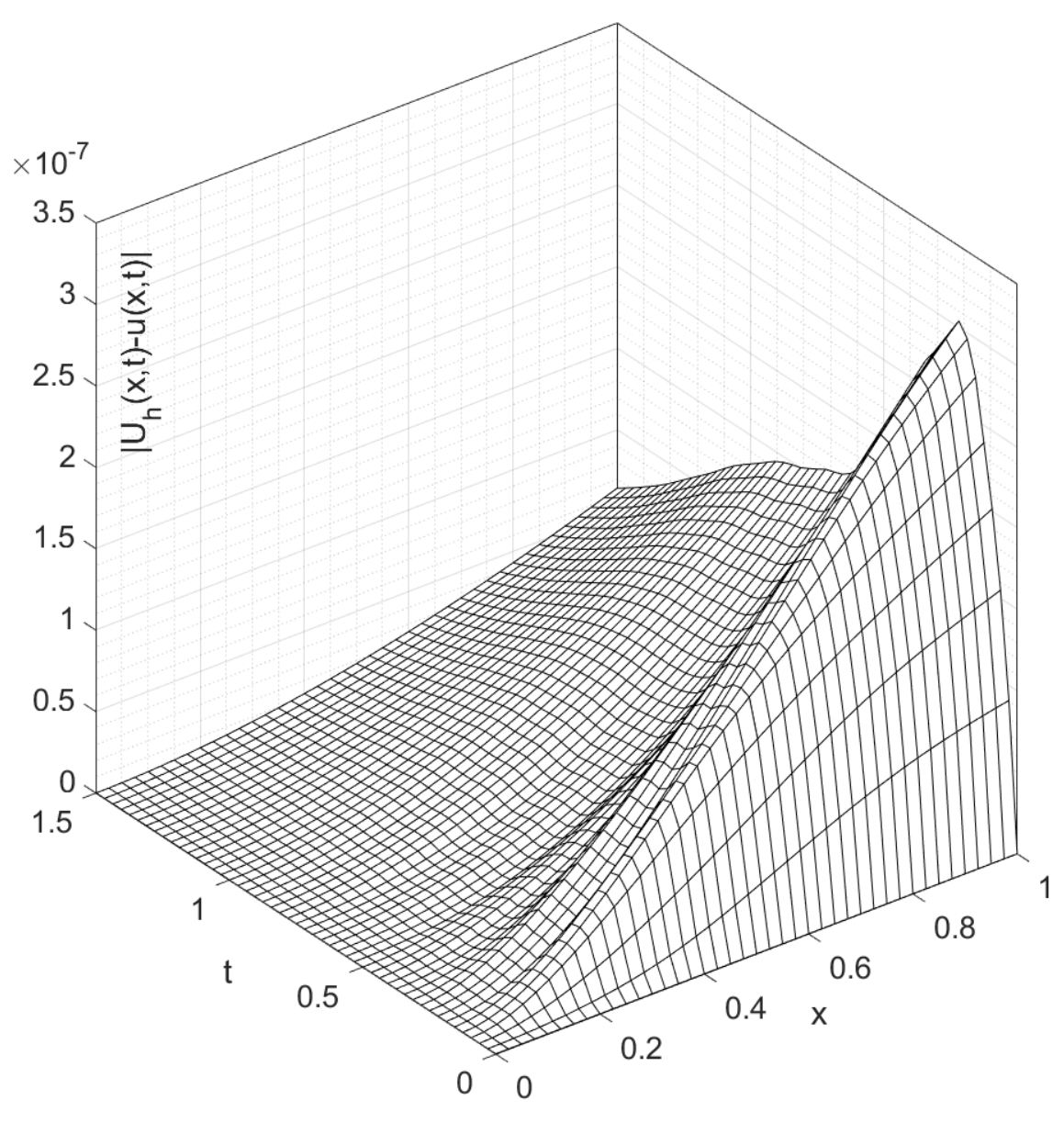}
	\caption{Numerical solution $U_h(x,t)$ and error  $|U_h(x,t)-u(x,t)|$.} \label{FigU}
\end{figure}

\begin{figure}[!htb]
	\includegraphics[width=6.8cm,height=6.8cm]{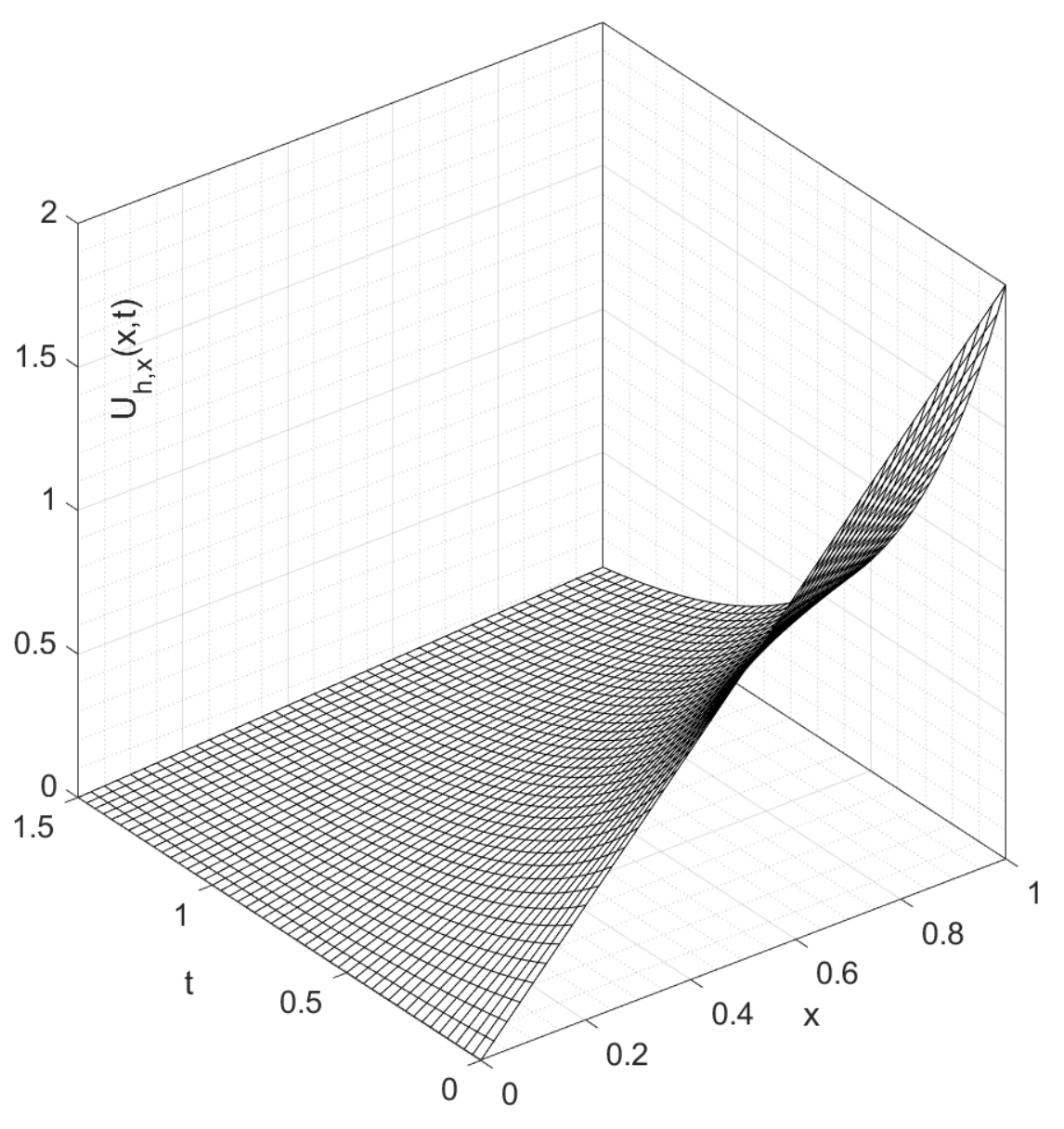}~~\includegraphics[width=6.9cm,height=6.9cm]{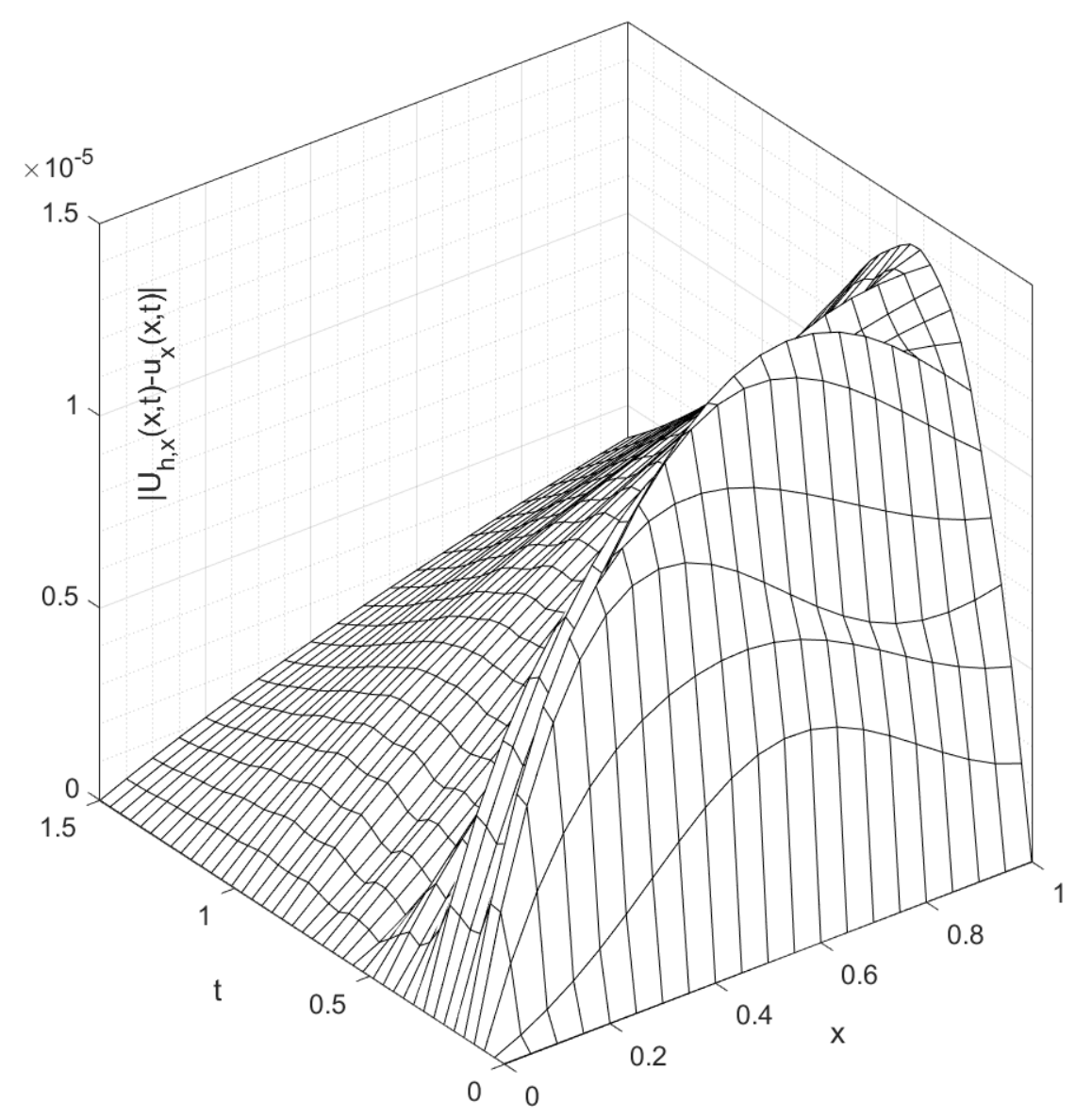}
	\caption{Numerical solution $U_{h,x}(x,t)$ and error  $|U_{h,x}(x,t)-u_x(x,t)|$.} \label{FigUx}
\end{figure}

The numerical studies below allow to analyze graphically the influence of the boundary control parameters on the stabilization of the beam vibration and on the asymptotic behaviour of the energy of the system. We also illustrate the verification of the theoretical results throughout the paper by this numerical test.
\begin{eqnarray}\label{NE1}
\left\{ \begin{array}{ll}
u_{tt}+2 u_{t}+\left((1+x)u_{xx}\right)_{xx}=0,\,~~ (x,t) \in \Omega_T:=(0,1)\times(0,1.5], \\ [4pt]
    u(x,0)=x^2, ~ u_t(x,0)=-2 x^2,\,~~ x\in (0,1),\\ [4pt]
	u(0,t)=u_{x}(0,t)=0,~~~t\in [0,1.5]\\
-(1+x)u_{xx}|_{x=1}= 6 u_x(1,t)+3 u_{xt}(1,t)=-4\exp(-2t),~~~t\in [0,1.5] \\ [4pt]
 \left((1+x)u_{xx}\right)_{x}|_{x=1}=4\, u(1,t)+2 u_t(1,t)=2\exp(-2t),~~~t\in [0,2],
\end{array} \right.
\end{eqnarray}
Here boundary spring parameters are $k_r=6,~~k_d=4$ and the damper parameters are $k_a=3,~~k_v=2$. The exact solution of (\ref{NE1}) and its first partial derivative with respect to $x$ are $\displaystyle u(x,t)=x^2\exp(-2t)$ and $\displaystyle u_x(x,t)=2\, x\, \exp(-2t)$. Numerical approximation of these functions can be found directly from the MOL technique in (\ref{fulldisc}) with the ratio of mesh parameters $h_x/h_t=40$. Corresponding approximate results are quite accurate  and illustrated in Fig. \ref{FigU} and Fig. \ref{FigUx}.

\begin{figure}[!htb]
	\includegraphics[width=6.8cm,height=6.8cm]{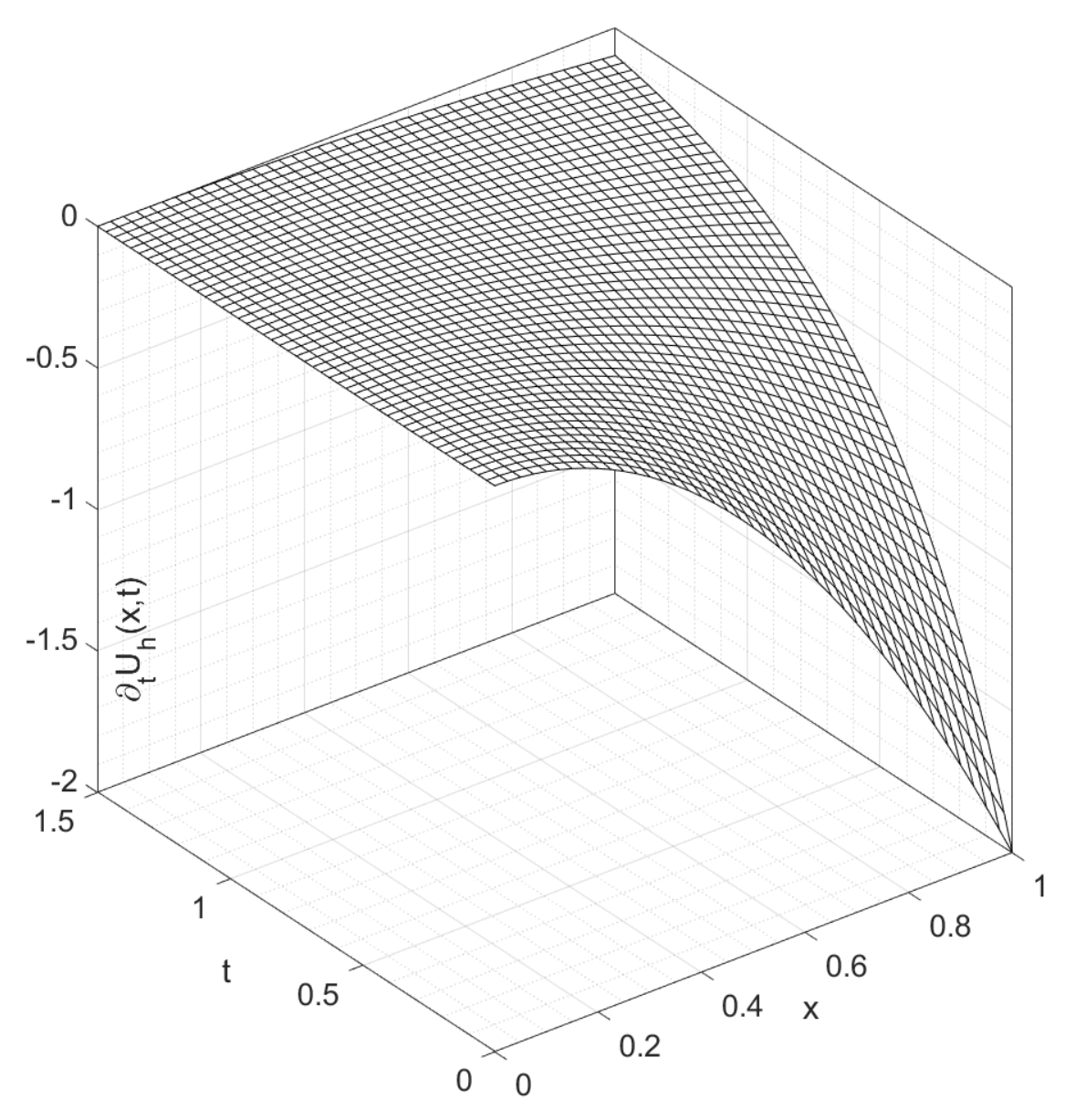}~~\includegraphics[width=6.9cm,height=6.9cm]{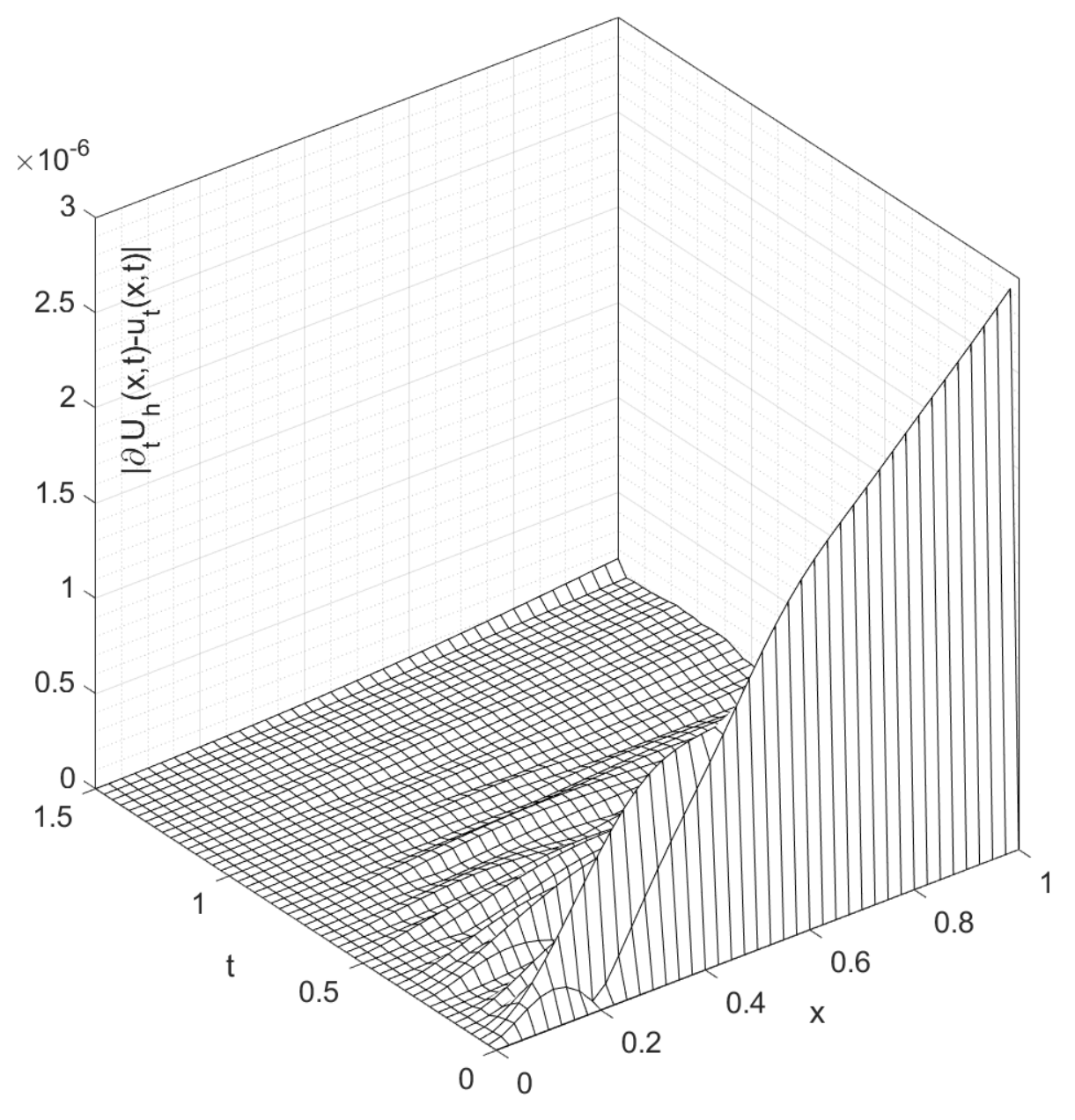}
	\caption{Numerical solution $\partial_tU_h(x,t)$ and error  $|\partial_tU_h(x,t)-u_t(x,t)|$.} \label{FigUt}
\end{figure}

\begin{figure}[!htb]
	\includegraphics[width=6.8cm,height=6.8cm]{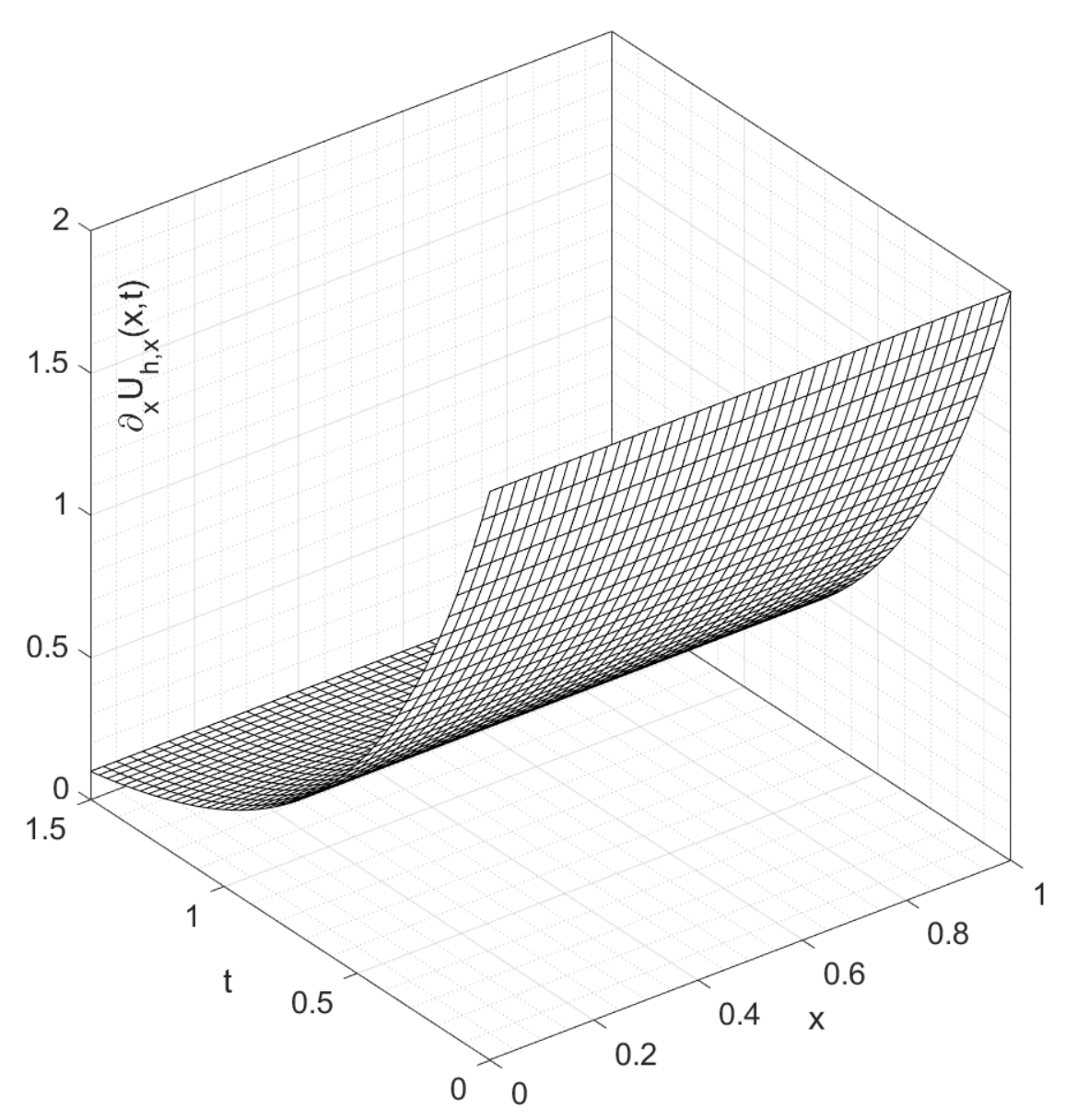}~~\includegraphics[width=6.8cm,height=6.8cm]{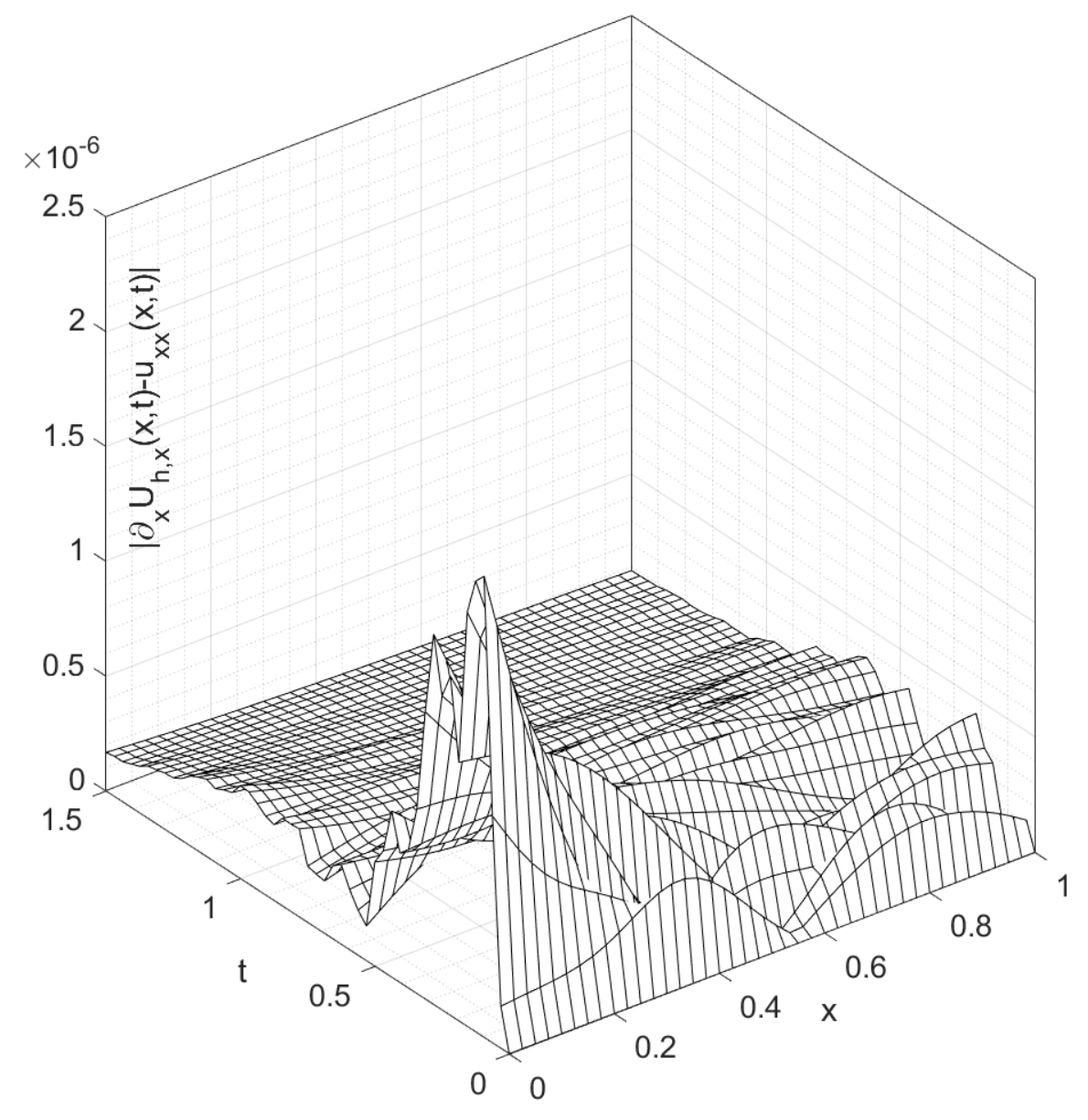}
	\caption{Numerical solution $\partial_xU_{h,x}(x,t)$ and error  $|\partial_xU_{h,x}(x,t)-u_{xx}(x,t)|$.} \label{FigUxx}
\end{figure}

\begin{figure}[!htb]
	\includegraphics[width=6.65cm,height=6.65cm]{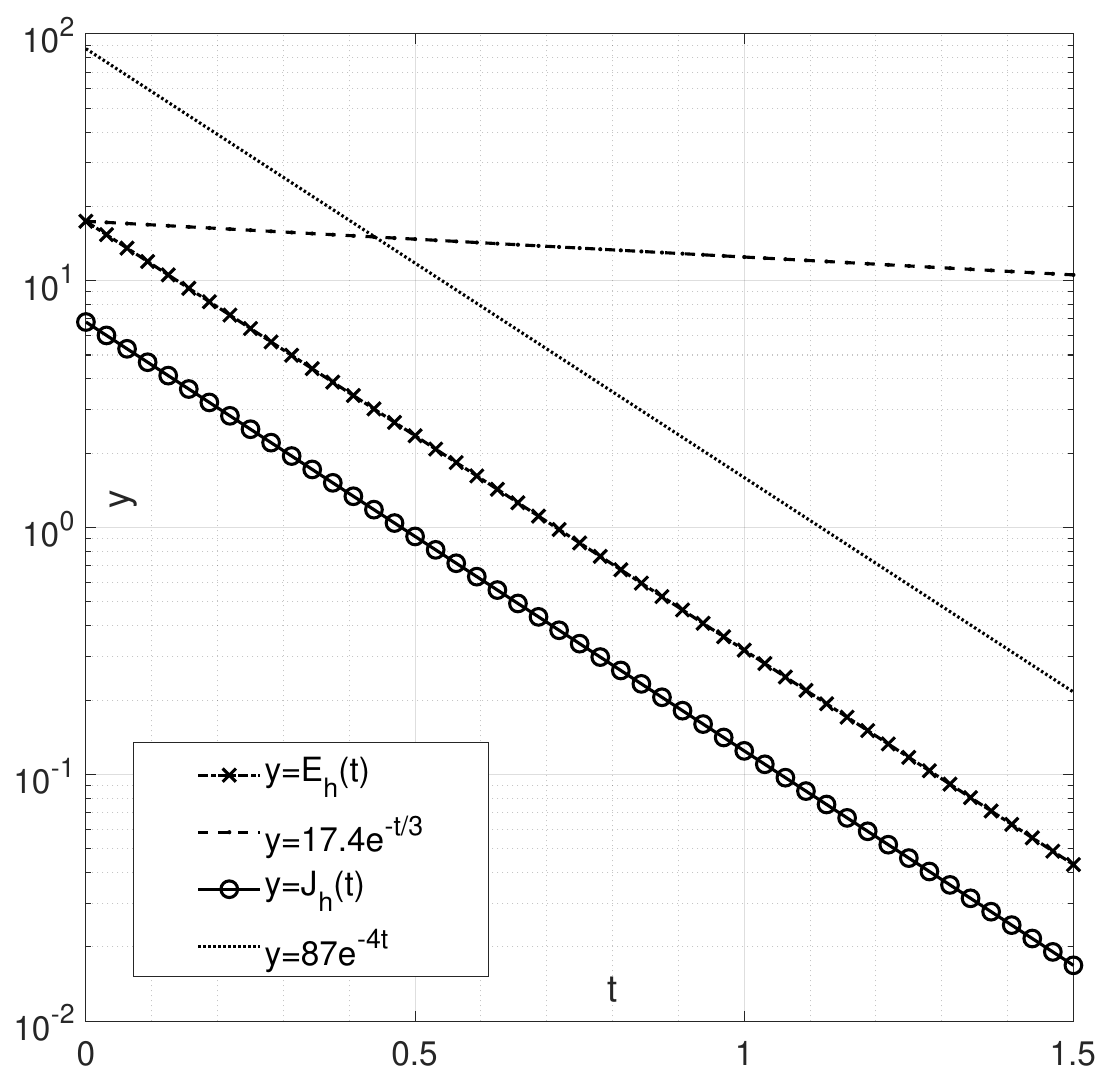}~~\includegraphics[width=6.9cm,height=6.9cm]{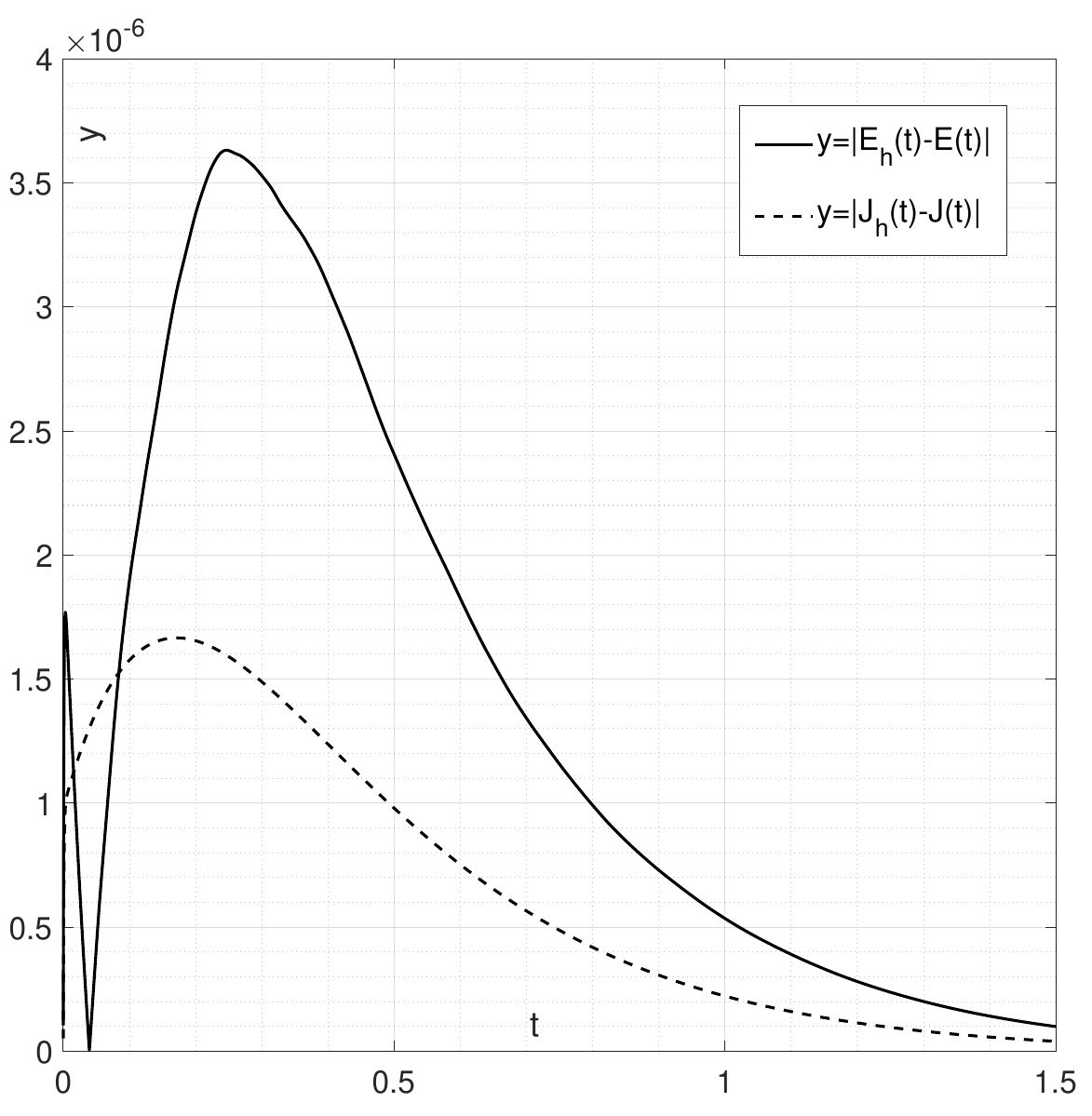}
	\caption{Approximation of  $\mathcal{E}_h(t)$ and $\mathcal{J}_h(t)$ with their upper bounds in logarithmic scale and errors $|\mathcal{E}_h(t)-\mathcal{E}(t)|$ and $|\mathcal{J}_h(t)-\mathcal{J}(t)|$.} \label{FigEJ}
\end{figure}

The energy function $\mathcal{E}(t)$ and auxiliary function $\mathcal{J}(t)$ for the given problem (\ref{NE1}) can be found respectively  $\displaystyle \mathcal{E}(t)=17.4 \exp(-4t)$ and  $\displaystyle \mathcal{J}(t)=6.8 \exp(-4t)$. In order to find their approximations $\mathcal{E}_h$ and $\mathcal{J}_h$  one needs to compute $u_t(x,t) $ and $u_{xx}(x,t)$. For this, we use centered difference quotient as follows.
\begin{eqnarray}\label{cent_dif}
\left\{ \begin{array}{ll}
\displaystyle u_t(x,t_j)\, \approx \, \partial_t\, U_h(x,t_{j}) =\left[U_h(x,t_{j+1})-U_h(x,t_{j-1})\right]/{2\,h_t},\\
\displaystyle u_{xx}(x_i,t)\, \approx \,\partial_{x}\, U_{h,x}(x_i,t) =\left[U_{h,x}(x_{i+1},t)-U_{h,x}(x_{i-1},t)\right]/{2\,h_x}.
\end{array} \right.
\end{eqnarray}
Approximate form of these derivatives $\partial_t\, U_h(x,t)$ and $\partial_{x}\, U_{h,x}(x,t)$ are obtained as a result of this centered difference approach and  shown in Fig. \ref{FigUt} and Fig. \ref{FigUxx} with their absolute errors.
Therefore, by replacing all of these approximations represented in Figs. \ref{FigU}-\ref{FigUxx} with corresponding exact quantities in (\ref{4}) and (\ref{11}), we obtain desired approximations $\mathcal{E}_h\, \approx \, \mathcal{E}(t)=17.4 \exp(-4t)$ and $\mathcal{J}_h\, \approx \, \mathcal{J}(t)=6.8 \exp(-4t)$. The accuracy of these approximations are illustrated in Fig. \ref{FigEJ} (right).\\

The upper bound of $\mathcal{J}(t)$ and  $\mathcal{E}(t)$ are follows from  (\ref{13}) and (\ref{22b}), respectively.
Here $\beta_0=1/2$ and $\beta_1=5$, then $$\mathcal{J}(t)=6.8 \exp(-4t)\leq \beta_1\, \mathcal{E}(t)=87 \exp(-4t).$$
Similarly, $\lambda=1$, $M_d\in(1,12)$ and $\sigma\in(0,1/3)$ for the considered test problem (\ref{NE1}). Therefore,
$$\mathcal{E}(t)=17.4 \exp(-4t) \leq 17.4 \exp(-t/3)< M_d\exp(-\sigma\, t) \mathcal{E}(0).$$
All these numerical studies related to $\mathcal{E}(t)$ and $\mathcal{J}(t)$  verify the theoretical results given in Proposition \ref{Proposition-3} and Theorem \ref{Theorem-1} and  are illustrated in Fig. \ref{FigEJ} (left).

\section{Some preliminary conclusions}

In this study we propose an approach which allows to obtain an explicit form of energy decay estimate for typical systems governed by Euler-Bernoulli beam controlled by boundary springs and dampers. As far as our knowledge extends, the relationship between the decay rate parameter $\sigma>0$ in the exponential stability estimate $\mathcal{E}(t)\le M_d\,  e^{-\sigma t}\, \mathcal{E}(0)$ and the physical parameters of the problem, including the damping parameters and the boundary dampers, was established here for the first time in the literature. This achievement was made possible through the utilization of a mathematical method rooted in the Lyapunov stability approach. It can be shown that in addition to the above studied cases, the considered approach is also applicable for cases of pinned-pinned, pinned-sliding, sliding-pinned, and sliding-sliding boundary conditions, including various types of inputs on the boundary $x=\ell$.


\end{document}